\documentclass[11pt,letterpaper]{amsart}
\usepackage{amsmath}%
\usepackage{amsfonts}%
\usepackage{amssymb}%
\usepackage{graphicx}
\usepackage{pictexwd,dcpic}
\usepackage{mathrsfs}
\usepackage{booktabs}
\usepackage{wrapfig}
\usepackage[margin=1in]{geometry}
\usepackage{float}
\usepackage [english]{babel}
\usepackage [autostyle, english = american]{csquotes}
\MakeOuterQuote{"}
\usepackage[lastpage,user]{zref}
\usepackage{fancyhdr}
\usepackage{wrapfig}

\title{A Gordian Pair of Links} 

\author{Rob Kusner}
\address{Dept. of Mathematics \& Statistics, University of Massachusetts, Amherst, MA 01003, USA}
\email{profkusner@gmail.com\\
kusner@math.umass.edu}

\author{W\"oden Kusner}
\address{
 Dept. of Mathematics, Vanderbilt University, Nashville, TN 37240, USA}
\email{w.kusner@vanderbilt.edu\\
wkusner@gmail.com}

\thanks{\noindent Work at the Aspen Center for Physics supported in part by National Science Foundation Award PHY-1607611; WK additionally supported by Austrian Science Fund (FWF) Project 5503 and NSF Award DMS-1516400.}

\date{}							
\begin{document}
\bibliographystyle{plain}
\maketitle
\vspace{-2em}
\begin{center}
{\scriptsize ABSTRACT: We construct a pair of isotopic link configurations that are not thick isotopic while preserving total length. }
\end{center}
\vspace{-.5em}
\noindent
\begin{figure}[h]
\centering
 \includegraphics[width=.46\textwidth]{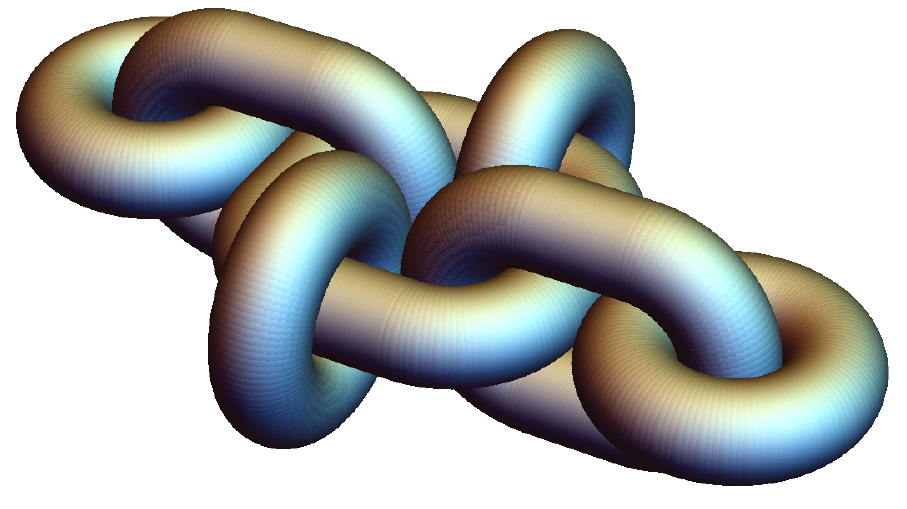}
 \centering
\includegraphics[width=.46\textwidth]{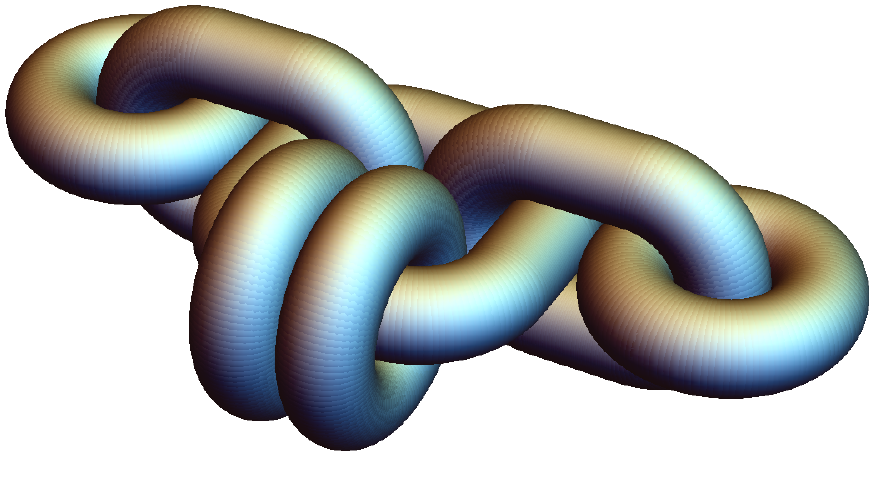}\\
{\scriptsize{\em A Gordian Pair:}  Configurations $\bf{R}$~(rotor) and $\bf{W}$~(wing) both minimize total ropelength in a common isotopy class; no isotopy between them preserves total ropelength.}
\end{figure}

Coward and Hass~\cite{coward2015topological}, using tools from~\cite{cantarella2002minimum}, gave an example of {\em physically distinct} isotopic configurations for a $2$-component link: 
No isotopy can be performed while preserving the ropelength of each component; 
however {\em length trading} among components, which is more natural in the criticality theory~\cite{cantarella2006criticality, cantarella2014ropelength} for ropelength, is not allowed.
Our configurations $\bf{R}$ and $\bf{W}$ are physically distinct in the broader length-trading sense appropriate for the Gordian unknot and unlink Problems~\cite{pieranski2001gordian}: {\em Do nontrivial ropelength-critical configurations of unknots and unlinks exist}?  
This problem arises in---and possibly obstructs---variational approaches~\cite{freedman1994mobius} to the Smale Conjecture~\cite{hatcher1983proof} via the space of unknots, and its generalization to spaces of unlinks~\cite{brendle2013configuration}.   

{\bf Definition.}  {\em A pair of link configurations is} Gordian {\em if the links are isotopic, but there is no isotopy between them with thickness at least $1$ which preserves total length.}

\noindent
In fact, we prove a stronger statement, in the context of {\em link homotopy} and {\em Gehring}~\cite{cantarella2006criticality} thickness:

{\bf Theorem.}  {\em The configurations $\bf{R}$ and $\bf{W}$ minimize total ropelength in their common link homotopy class, but there is no link homotopy between them with Gehring thickness at least $1$ while preserving the total ropelength.}
\newline
\noindent
Because link homotopy is coarser than isotopy, and since the Gehring thickness constraint is more permissive than that for standard~\cite{cantarella2002minimum} thickness, {\em a fortiori} this is a Gordian pair.

{\it Proof of Theorem.} (i) For any minimizing configuration in this link homotopy class, each component must be a particular type of stadium curve~\cite{cantarella2006criticality, cantarella2002minimum} surrounding 1, 2 or 4 disjoint unit disks; in the last case, there is an interval moduli space of such curves $C$, ranging between the square (depicted above) and equilateral-rhombic configurations of 4 unit disks.
(ii) Define a map $\Pi$ from the space of minimizing link configurations in this link homotopy class to the space $\mathcal{C}_4(S^1)$ of 4-point configurations on the circle, taking the given link configuration to the 4 intersection points of $C$ with the planar spanning disks for the 4 components linking $C$.
(iii) The image of $\Pi$ lies in the closed subset of $\mathcal{C}_4(S^1)$ where each intersection point lies in one of the 4 curved arcs of $C$, a deformation retract of $\mathcal{C}_4(S^1)$.
(iv) The 4-configurations  $\Pi(\bf{R})$ and $\Pi(\bf{W})$ lie in distinct path components of $\mathcal{C}_4(S^1)/\textrm{O}(2)$---corresponding to dihedral orders of 4 points on a circle---so there is no path between $\bf{R}$ and $\bf{W}$ in the moduli space of Gehring-ropelength minimizers.  \qed

{\bf Remark.} In forthcoming work (in part with Greg Buck), we develop tools giving a stronger result: {\em The total Gehring ropelength must rise by at least~$2$ in any isotopy (or link homotopy) between these minimizing link configurations.} 

\bibliography{rope}
\end{document}